# Reviving networked multi-dimensional dynamical systems

Nan Dong[1,2,3], An Zeng[3*], Honggang Li[1,2,3**]

[1]Department of Systems Science, Faculty of Arts and Sciences, Beijing Normal University, Zhuhai 519087, China

[2]International Academic Center of Complex Systems, Beijing Normal University, Zhuhai 519087, China

[3]School of Systems Science, Beijing Normal University, Beijing 100875, China

e-mail: anzeng@bnu.edu.cn; hli@bnu.edu.cn

**Abstract**: From gene regulatory networks to mutualistic networks, controlling a single node in the network topology can transform these complex dynamical systems from undesirable states to desirable ones. Corresponding methods have been well-studied in one-dimensional dynamical systems. However, many practical dynamical systems require description by multi-dimensional dynamical systems, such as the mutualistic symbiotic systems formed by flowering plants and pollinating insects. Existing one-dimensional methods cannot handle the cases of multi-dimensional dynamical systems. Based on this, we propose a method to control a single node to activate network connections in multi-dimensional dynamical systems. In such systems, the changes of each node are described by multiple nonlinear differential equations. All remaining nodes are stratified according to the shortest path to the controlled node, thereby reducing the dimensionality of the system. Such a large-scale dynamical system can ultimately be replaced by a very simple system. By analyzing the reduced-dimensional system, we can predict the extent of control needed to restore the system state. We apply this method to a wide range of fields, achieving activation of various real multi-dimensional complex dynamical systems.

# 1 Model Framework

## 1.1 Introduction

Complex systems often have an ideal state in which they function properly and a non-ideal state where they do not[1]. Examples include ecosystems[2], the regulation and expression of genes and proteins within cells[1], financial markets[3], and labor markets. These systems commonly face the challenge of how to recover from a non-ideal state to the ideal one. This paper refers to the process of a complex system recovering from a non-ideal state to an ideal state as activation. One obvious approach to activation is to control all parts of the system. However, since real-world complex systems are often very large, this method can be costly or even unfeasible. Another approach is to control a part of the system and, by utilizing the interconnections between different parts, activate the whole system. Many examples are similar to the latter approach. For instance, it is common for the same theory or invention to be proposed multiple times. The topology and level of support of a network can influence or even determine whether it remains unnoticed or becomes widely known. Examples of such phenomena include the theory of evolution, chaotic systems, and the scale-free properties of networks.

It is necessary to predict, using quantitative methods, whether local control can activate the system. Such a system can be abstracted as a networked dynamical system. The changes of each node are represented by differential equations. However, there are significant challenges in dealing with this. Real networks are often very large in scale, their topological structures are difficult to determine, and the changes in each node may require multidimensional equations to be represented. Often, we can only know the macroscopic properties such as the scale of the network and its average degree. Therefore, an ideal approach would be one that only requires knowledge of the network's characteristics and has low computational complexity.

In 2022, Hillel et al. addressed the problem of networked one-dimensional dynamical systems under the assumption of an infinitely large network[1]. They conducted extensive numerical simulations to validate the effectiveness of their method, which also proved to be valid for real-world network topologies. However, their approach has some limitations. For instance, it requires the network to be locally tree-like, the dynamical equations to be separable, and the control to be constant. Moreover, the time it takes for the system to reach a steady state is unknown. The main limitation of this method is that it cannot provide a general solution for multi-dimensional dynamical systems. To describe the changes in complex systems, it is sometimes necessary to use multi-dimensional dynamical systems that capture interactions between different groups. For example, gene expression and regulation involve DNA, RNA, and proteins; mutualistic ecosystems include flowering plants and pollinating insects; and labor market networks involve unemployment and employment.

Based on this, we propose a method to solve the problem of predicting whether a microscopic intervention can activate a networked multi-dimensional dynamical system. The reason why the method proposed by Hillel et al[1]. cannot handle multi-dimensional cases is that, when applied to multi-dimensional systems, it results in a set of nonlinear multivariate equations, which do not have a general analytical solution. Numerical methods are also very challenging. Therefore, we take a different approach, assuming a finite number of nodes in the network. Apart from the selected node, the remaining nodes are grouped into layers based on the shortest path to the selected node, thus reducing the dimensionality of the system. The connections between nodes in different layers are then inferred from the network's scale and topological information. In real networks, the maximum shortest path between any two nodes is usually proportional to the logarithm of the number of nodes. For example, it only takes at most six edges to connect any two people in the world. This suggests that even for very large-scale systems, they can still be simplified into a simple model.

1.2 Dynamic Systems

The dynamic system we are dealing with can be represented by the following equation (1.1):

$$\frac{dX_i(t)}{dt} = G_i(X_i(t)) + F_i[X_1(t), X_2(t)...X_N(t); A]. \quad (1.1)$$

Here, $A$ represents the adjacency matrix with a size of $N \times N$. $A_{ij} = 1$ indicates that there is an edge between two nodes, while $A_{ij} = 0$ indicates that there is no edge between them, where $i, j = 1, 2 ... N$. This is an undirected network, so the matrix $A$ is symmetric. $X_i(t)$ represents the state of node $i$ at time $t$. It is a vector and can represent multiple activities. $G_i(X_i(t))$ represents the effect of the node's own activities on its change, while $F_i[X_1(t), X_2(t)...X_N(t); A]$ represents the effect of other nodes connected to node $i$ on its change. Equation (1.1) is not intuitive enough, so we will use a specific example to describe how to handle such problems.

$$\frac{dx_i(t)}{dt} = -B_1 x_i(t) + \sum_{j=i}^{N} A_{ij} \frac{y_j^2(t)}{1 + y_j^2(t)},$$
$$\frac{dy_i(t)}{dt} = -B_2 y_i(t) + \sum_{j=i}^{N} A_{ij} \frac{x_j^2(t)}{1 + x_j^2(t)}. \quad (1.2)$$

Here, both parameters $B_1$ and $B_2$ are positive numbers. This is a model that describes gene regulation on a complex network. We will provide a detailed introduction to this model later.

1.3 Network Topology

Suppose the total number of nodes in the network $N$, and the average degree $k$, are known, with $k \ll N$. The reason we only need these two parameters is based on practical considerations. In real-world networks, such as the human brain network, the connection matrix is difficult to obtain, and the specific values are constantly changing. By using only these two values, the cost of obtaining information is kept low. The network is undirected. In real scenarios, some networks may be directed or weighted, and these cases can also be addressed within this framework. In the case of weighted networks, the expected value can be used as a substitute. In the case of directed networks, it is necessary to modify the section in Section 3 that derives the average number of edges per layer of the network.

The network is connected. If the network is not connected, the selected node cannot influence parts of the network that are not connected to it. Between any two nodes in the network, there is at most one directly connected edge, meaning there are no cycles. In an undirected network, the shortest path length between any two nodes $i$ and $j$ is defined as the distance $L_{ij}$ between them. For example, if $L_{ij} = 2$, it means that at least two edges are required to travel from node $j$ to node $i$. In the complex network, we randomly select a node for control, denoted as $s$. Using the distance definition $L_{ij}$, we layer all nodes, except for node $s$, based on their distance to node $s$.

$$K_s(l) = \{i = (1,...,N) | L_{si} = l\}. \quad (1.3)$$

Equation (1.3) defines the set of all nodes at a distance $l$ from node $s$. For example $K_s(0) = \{s\}$,

$K_s(1)$ represents the sets of all nodes directly connected to node $s$ by an edge. Based on this definition, we can rewrite equation (1.2) as follows:

$$\frac{dx_i(t)}{dt} = -B_1 x_i(t) + \sum_{j \in K_i(1)} \frac{y_j^2(t)}{1 + y_j^2(t)},$$
$$\frac{dy_i(t)}{dt} = -B_2 y_i(t) + \sum_{j \in K_i(1)} \frac{x_j^2(t)}{1 + x_j^2(t)}. \quad (1.4)$$

Here, $j \in K_i(1)$ represents the set of all nodes directly connected to node $i$.

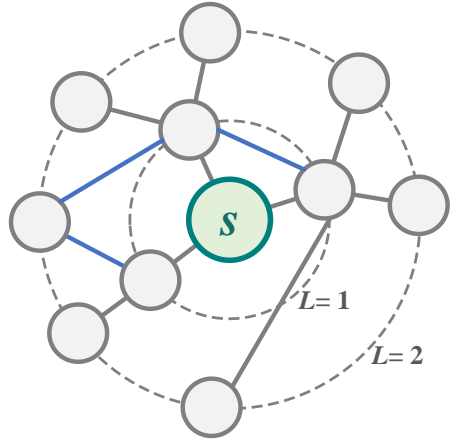

Fig.1.1 shows the layering of all nodes in the network except for the selected node $s$. The layer of a node is determined by its shortest path to the selected node $s$. Nodes within the same layer can be connected by edges, as shown by the blue edges in the first layer. Outer layer nodes can have two edges from inner layer nodes, as illustrated by the nodes in the second layer.

In order for the analysis to proceed, the entire dynamical system needs to satisfy certain conditions. Some of these conditions are strict, and not meeting them may lead to a decrease in prediction accuracy. Others are less strict, intended only for convenience in later derivations.

1 When reaching a steady state, the activity of each node in the system should be narrowly distributed. This is not a requirement specific to the network or the dynamical model alone, but rather a combined requirement for both. A typical case that does not meet this condition is when the node activities follow a power-law distribution. This is usually caused by the power-law distribution of degrees in the network, which can lead to certain phenomena in specific dynamical models. However, a power-law degree distribution does not necessarily result in node activities following a power-law distribution. This condition is mandatory, and failure to meet it will decrease the accuracy of predictions. The reasons for this will be explained in the numerical simulation section later.

2 The topology of different nodes should not have essential differences. That is to say, when selecting a node at random multiple times, the connection characteristics of nodes at the same level are similar. Networks with central-peripheral structures or community structures do not meet this condition. This condition ensures that, in the dimensionality reduction process described later, the nodes within the same layer can be replaced by an average node. Networks that do not meet this condition can sometimes still be processed according to the approach presented here, but the connections in the network would need to be derived again.

3 The connections between nodes are random. Furthermore, we assume that there is at most one edge between any two nodes in the network. That is to say, the network has no loops. Under the condition of random connections, the probability of this occurrence is particularly low. If the probability of this occurring is relatively high, the network can be treated as a weighted network. This assumption is made for the purpose of deriving the average number of edges between different layers later. This assumption does not need to be strictly satisfied. If not, it can be addressed using statistical methods or shortest-path distribution techniques.

1.4 System State

Real complex systems typically have multiple steady states, and the corresponding dynamical systems of complex network connections also have multiple equilibrium points. Without loss of generality, we denote the multiple equilibrium points obtained from system (1.2) as $X^\alpha = (x_{1,\alpha}, y_{1,\alpha}; x_{2,\alpha}, y_{2,\alpha}; \ldots; x_{N,\alpha}, y_{N,\alpha})$, where $\alpha = 0, 1, 2, \ldots$ represents the number of equilibrium points. Steady states can be determined by setting the right-hand side of equation (1.2) equal to 0. We assume that the values of non-ideal states are either 0 or tend toward 0. The remaining equilibrium points represent non-zero states. This assumption is crucial for our theory. This point has not been adequately addressed in the research by Hillel et al[1]. In fact, without such a condition, the method proposed by Hillel et al. would fail.

The reason we make this assumption is based on the following considerations: 1 The study of the transition from a non-zero equilibrium state to the zero equilibrium state has been explored by others, with a focus on the network topology perspective. 2 In practical situations, it is rare for a non-zero equilibrium state to transition to the zero equilibrium state through microscopic interventions. For example, widespread infectious diseases rarely disappear through the control of only a few individuals. 3 When the entire system is in a non-zero equilibrium state, each node's influence comes not only directly or indirectly from the selected nodes but also from other nodes, which complicates the analysis. In summary, although many phenomena may arise from small factors, these phenomena usually remain stable. Only when the system structure undergoes significant changes may these phenomena disappear. Therefore, we do not analyze the issue of microscopic interventions causing the entire system to transition from a non-zero equilibrium state to the zero equilibrium state. Additionally, we consider only the case where the system is in a single steady state. Currently, some researchers have observed the coexistence of multiple steady states in both theoretical models and real systems. For example, in the East African savannah, multiple steady states involving different grass and tree species can coexist under similar environmental conditions. Since the mechanisms of such phenomena are not yet fully understood, we do not consider systems where this occurs.

The degree of a node has a crucial impact on the activity of the node. Since the network's adjacency matrix is heterogeneous, the activity of each node is also quite complex. However, in the case of random connections, most nodes' degrees are concentrated around the mean degree. Complex systems often have multiple steady states. In the absence of external interventions, non-ideal steady states typically do not spontaneously transition to ideal steady states. Therefore, our goal is to use microscopic interventions to shift the system from a non-ideal steady state to an ideal steady state, or to bring the system into the attractor domain of the ideal steady state. Without loss of generality, we denote the non-ideal steady state as $X^0$ and the ideal steady state as $X^1$. Thus, our objective is to find the appropriate control such that, after controlling the system for a sufficiently long time, it will enter the attractor domain of the ideal steady state.

$$\mathfrak{R}_1 = \left\{ X^0(t=0) \mid X^1(t \to \infty) \right\}. \quad (1.5)$$

Here, $\mathfrak{R}_1$ contains all the control combinations that, after being applied for a sufficiently long time, can shift the system's equilibrium state from $X^0$ to $X^1$. Our goal is to determine, given the network topology and the corresponding dynamical system, which controls can drive the system to the ideal equilibrium state.

2 Activation of the System

2.1 Nodes for Microscopic Intervention

Assume that a dynamical system described by equation (1.1) has two stable states: a non-ideal steady state and an ideal steady state. We now intervene microscopically to change the system from the non-ideal steady state to the ideal steady state. During the activation process, because only a few nodes are selected for control, the specific selection method can be divided into three types. The first type is to select one node. In this case, the system will either activate or the control effect is almost limited to the controlled node, with little change in the states of the remaining nodes. The second type is to randomly select a few nodes for control. Due to the large scale of the network, these nodes are usually far apart and not directly connected. The effect of the controlled nodes is isolated and cannot influence each other. This situation is equivalent to the first one. The third type is to select a few adjacent nodes for control. In this case, more nodes directly connected to the selected ones are affected, making the situation more complex. The control effects on different nodes can accumulate. The difference between the third and first types is that the number of edges directly connected to the controlled nodes changes. This situation can be addressed using this framework. When we apply control, the state change of the system is shown in equation (2.1).

$$\frac{dX_s(t)}{dt} = \delta_s(t) \qquad\qquad\qquad s \in \mathfrak{R},$$
$$\frac{dX_i(t)}{dt} = G_i(X_i(t)) + F_i[X_1(t), X_2(t)...X_N(t); A] \quad i \notin \mathfrak{R}. \qquad (2.1)$$

$\delta_s(t)$ represents the control function of the controlled nodes. In addition to the selection of nodes, another issue is the method of control. During the entire control process, the control applied to the selected nodes can be varied in different ways, such as fixed control or varying control. This method can address the case of multiple nodes with varying control, which is not covered by the method of S et al.

2.2 Dimensionality Reduction of the Network

After clarifying the necessary conditions and explanations, we proceed with the dimensionality reduction of the entire system. This reduction process reflects the fact that, during control, the activities of nodes at the same layer may differ, but they follow the same distribution and can be replaced by an averaged node. This is easy to understand intuitively. Since there is no significant difference in the topology of each node, the state of each node depends on its distance from the selected node. Thus, after layering the system, it can be represented as:

$$\frac{dX_s(l,t)}{dt} = \frac{1}{|K_s(l)|} \sum_{i \in K_s(l)} \frac{dX_i}{dt}. \quad (2.2)$$

Here, the symbol '| |' represents the number of elements in a set, so '$|K_s(l)|$' denotes the total number of nodes whose distance from the selected node $s$ is $l$. '$X_s(l,t)$' represents the average state of all nodes whose distance from the selected node $s$ is $l$. Therefore, equation (2.1) can be further expressed as:

$$\frac{dX_s(l,t)}{dt} = \frac{1}{|K_s(l)|} \sum_{i \in K_s(l)} \left( G_i(X_i(t)) + F_i[X_1(t), X_2(t)...X_N(t); A] \right). \quad (2.3)$$

Using (2.3), equation (1.4) can be expressed as:

$$\frac{dx_s(l,t)}{dt} = -B_1 x_s(l,t) + \frac{1}{|K_s(l)|} \sum_{i \in K_s(l)} \sum_{j \in K_i(1)} \frac{y_j^2(t)}{1+y_j^2(t)},$$
$$\frac{dy_s(l,t)}{dt} = -B_2 y_s(l,t) + \frac{1}{|K_s(l)|} \sum_{i \in K_s(l)} \sum_{j \in K_i(1)} \frac{x_j^2(t)}{1+x_j^2(t)}. \quad (2.4)$$

Here, since $x_i(t)$ and $y_i(t)$ are linear, $\frac{1}{|K_s(l)|}\sum_{i \in K_s(l)} B_1 x_i(t) = B_1 x_s(l,t)$ and $\frac{1}{|K_s(l)|}\sum_{i \in K_s(l)} B_2 y_i(t) = B_2 y_s(l,t)$ are obviously valid. As mentioned earlier, we require that the activity of each node does not fluctuate significantly, and the activity of nodes at each layer follows the same distribution. Therefore, for the $l$-th layer, there is $\frac{1}{|K_s(l)|}\sum_{i \in K_s(l)} \frac{y_i^2(t)}{1+y_i^2(t)} \approx \frac{y_s^2(l,t)}{1+y_s^2(l,t)}$. Many operations are similar and will not be repeated. As previously mentioned, in an undirected network, nodes connected to the l-th layer can only come from the $l-1$, $l$, or $l+1$ layers. Thus, equation (2.5) can be expressed as:

$$\frac{dx_s(l,t)}{dt} = -B_1 x_s(l,t) + \sum_{n=-1}^{1} c_{l,n} \frac{y_s^2(l+n,t)}{1+y_s^2(l+n,t)},$$
$$\frac{dy_s(l,t)}{dt} = -B_2 y_s(l,t) + \sum_{n=-1}^{1} c_{l,n} \frac{x_s^2(l+n,t)}{1+x_s^2(l+n,t)}.$$
(2.6)

Here, the parameter $c_{l,n}$ represents the average number of edges connecting each node in the $l$-th layer to nodes in the $l+n$-th layer. Clearly, $c_{l,-1} + c_{l,0} + c_{l,1} = k$. Equation (2.6) shows that, under certain conditions, a complex networked dynamical system can be represented by a simplified system that is proportional to the number of layers. Real networks and most model networks exhibit a characteristic where the maximum shortest path in the network is proportional to the logarithm of the number of nodes. Therefore, in most cases, these complex systems can be replaced by relatively simple systems.

3. The connection relationships between nodes at different levels

In the simplified model (2.6), only the connection parameter $c_{l,n}$ is unknown. For networks that satisfy the previous requirements, we will give the method to calculate $c_{l,n}$ below. This section is mainly focused on analyzing the case of selecting only one node. The case of selecting a few adjacent nodes can be solved using a similar process. The system's dimensionality reduction process replaces the degree of each node with the average degree $k$. The degree of most nodes in a complex network is close to the average value. Even for nodes with a degree deviating from the average, the main error in simplification comes from the outer layer, and the change in node state is mainly influenced by the states of the inner layer nodes. Therefore, this approach is reasonable.

To calculate $c_{l,n}$, we define the following parameters: the total number of nodes in layer $l$ is $d_l$, the total number of edges from layer $l-1$ to layer $l$ is $e_l$, the total number of remaining nodes outside layer $l$ is $f_l$, and the number of nodes in layer $l$ that have two edges from layer $l-1$ is $g_l$. These parameters are set to assist in the calculation of $c_{l,n}$.

A node is randomly selected and denoted as $s$. All remaining nodes are stratified based on their shortest path to node $s$. If the shortest path of a node is $l$, then the node resides in the $l$-th layer. Node $s$ has $k$ edges, hence the first layer contains $d_1 = k$ nodes. For each node in the first layer, there is one edge from the selected node $s$, so each node must distribute its remaining $k-1$ edges to the other $k-1$ nodes within the same layer or to the $f_1 = N - k - 1$ nodes that have not yet been connected. Since the connections are random, each node in the first layer has $c_{1,0} = \frac{(k-1)(d_1-1)}{(f_1+d_1-1)} = \frac{(k-1)^2}{N-2}$ edges connected to other nodes in the first layer and $c_{1,1} = k - c_{1,-1} - c_{1,0} = \frac{(k-1)(N-k-1)}{N-2}$ edges connected to nodes in the second layer.

Since the first layer has $d_1 = k$ nodes, the total number of edges between the first layer and the second layer is $e_2 = c_{1,1}d_1 = \frac{k(k-1)(N-k-1)}{N-2}$. Now, we analyze how many nodes in the second layer have two edges connected to the first layer. When connecting the first layer to the remaining unconnected nodes, there are $f_1 = N - k - 1$ unconnected nodes in total. Therefore, for each unconnected node, the probability of connecting to the first layer is $\frac{c_{1,1}d_1}{f_1}$, so the probability of each node having two edges is $\frac{(c_{1,1}d_1)^2}{f_1^2}$. Thus, the number of nodes in the second layer with two edges is $g_2 = \frac{(c_{1,1}d_1)^2}{f_1}$. The number of nodes with three edges, four edges, or more is very small and can be ignored. The number of nodes in the second layer is the total number of edges sent outward from the first layer minus $g_2$, which is $d_2 = c_{1,1}d_1 - g_2$. The remaining number of nodes in the second layer is $f_2 = f_1 - d_2$. The average number of edges between the second layer and the first layer is $c_{2,-1} = \frac{e_2}{d_2}$. The average number of edges within the second layer is $c_{2,0} = \frac{(k - c_{2,-1})(d_2 - 1)}{(f_2 + d_2 - 1)}$. We express the general formulas for the parameters of the $l$-th layer in a similar process: $e_l = c_{l-1,1}d_{l-1}$, $g_l = \frac{(c_{l-1,1}d_{l-1})^2}{f_{l-1}}$, $d_l = c_{l-1,1}d_{l-1} - g_l$, $f_l = f_{l-1} - d_l$, $c_{l,-1} = \frac{e_l}{d_l}$, $c_{l,0} = \frac{(k - c_{l,-1})(d_l - 1)}{(f_l + d_l - 1)}$, $c_{l,1} = k - c_{l,-1} - c_{l,0}$. Some parameter expressions contain other parameters from the same layer, so the parameters need to be calculated sequentially in the order given above.

In networks with a limited scale, the number of layers will also be limited. However, the above parameters can continue to be solved indefinitely, which is obviously not possible. Therefore, we need to set a termination condition. When the parameter of the $l$-th layer becomes negative, it indicates that the $l$-th layer is the outermost layer, and we should calculate all parameters according to the termination layer formula. At this point, $c_{l,-1} = \frac{c_{l-1,1}d_{l-1}}{f_{l-1}}$, $c_{l,0} = k - c_{l,-1}$. Since we have reached the outermost layer, $c_{l,1} = 0$. No further parameters need to be computed. In this process, besides the average number of edges per node at each layer, we also obtain the number of layers of the network. Due to the randomness of the network, the true number of layers may not exactly match the predicted number of layers. However, the impact of this on the prediction result can be ignored.

Through the above process, for a given dynamic system, we have obtained a simplified system. Even if the original system is very large, the simplified system is only proportional to the logarithm of the total number of nodes in the original system. The above method can be extended to weighted directed networks. For networks with more complex connection rules, the derivation process will be more difficult. This method is still applicable to networks that do not strictly meet the requirements. For networks or model networks that significantly deviate from the assumptions in the text, sampling, statistics, or methods such as shortest path distributions can be used to obtain these parameters.

4 Numerical Simulation

In this section, we validate the effectiveness of the method by using improved versions of several existing models. The reason for not directly using existing models is that higher levels of control are more commonly required to activate the system. However, we tested many existing dynamic models and found that they are typically either not activated by any control or are activated by any control. In these cases, the prediction accuracy of the method is very high, making it difficult to demonstrate the predictive effects. This also suggests that there are many unreasonable aspects in the existing models.

4.1 Explanation of Some Aspects of the Numerical Simulation

Before the numerical simulation, it is important to explain the reasons behind some of the settings in this section. The first key question is how to judge whether activation is successful. By definition, if the control allows the system to enter the basin of attraction of the ideal state or to reach the ideal state, the system is considered activated. However, determining the basin of attraction of the ideal state is not easy. Therefore, after applying the control, all nodes should evolve according to their own dynamics for a sufficiently long time to reach a stable state before making the judgment. During this second period, the changes in each node are no longer solely influenced by the selected node, but also by the other nodes. In this case, this method cannot predict the changes during this time.

Alternatively, if we use the average state of the outermost nodes at the end of the control as the basis for judgment, although their state is only minimally influenced by the controlled nodes, we would still need to identify which nodes belong to the outermost layer in the code, which complicates the calculation. Therefore, we judge the success of activation based on the average activity of all nodes at the end of the control. Despite the presence of controlled nodes, this is not exactly the same as the equilibrium state, but the difference is very small.

In addition, another point is the determination of the simulation time. The system can only reach equilibrium if the control lasts long enough. However, long numerical simulations consume significant computational resources. Therefore, the simulation time used below is the time at which the system just reaches equilibrium.

4.2 Gene Regulation Model

Equation (4.1) describes the regulatory system of RNA and proteins within a cell. In this process, there is positive feedback regulation between RNA and proteins. Many RNA and protein regulatory interactions within the cell can be represented by the following model. Existing studies show that the amount of protein in the cell is crucial for maintaining the normal function of the organism.

$$\frac{du_i(t)}{dt} = -B_1 u_i(t) + \sum_{j=i}^{N} A_{ij} \frac{v_j^2(t)}{1+v_j^2(t)},$$
$$\frac{dv_i(t)}{dt} = -B_2 v_i(t) + \sum_{j=i}^{N} A_{ij} \frac{u_j^2(t)}{1+u_j^2(t)}.$$
(4.1)

In the equation, $u_i(t)$ and $v_i(t)$ represent the amounts of RNA and protein at node $i$ at time $t$, respectively. B1 and B2 represent their degradation rates, with the degradation rate being proportional to their respective amounts. There is a nonlinear positive feedback regulation between RNA and protein. When the levels of RNA or protein are low, the regulatory effect is weak. As the amounts increase, the regulatory effect significantly strengthens. However, the regulatory effect of their increasing amounts does not rise indefinitely but has an upper bound.

The parameter values are set as $B_1 = 1.3$ and $B_2 = 1.5$. The model is verified on a network with 5000 nodes and an average degree of 10. Fixed control is applied, meaning that the control over the selected nodes remains unchanged throughout the process. The control time is set to 60. The results are shown in Figure 4.1. Each parameter combination is simulated 10 times, and the proportion of successful activations is recorded. The generated networks and selected nodes are random each time. It is important to note that, for this model, the node activity in the BA network follows a power-law distribution, which does not meet our requirements. The power-law distribution of node activity is a common phenomenon in real systems, and we aim to observe the performance of this method under conditions where the strict requirements are not met.

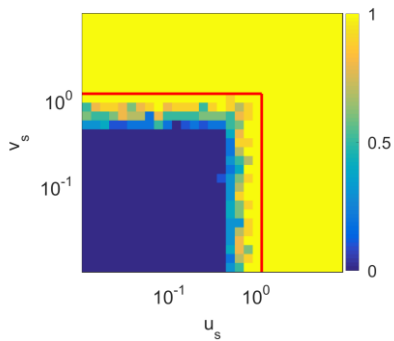 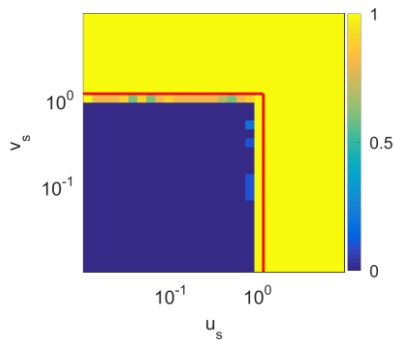 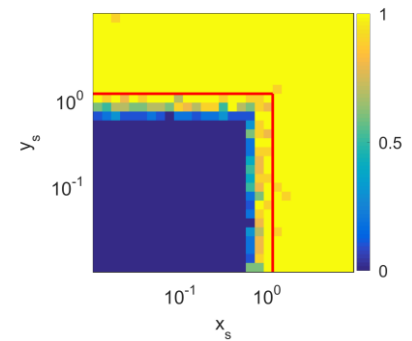

Fig. 4.1. (a) Prediction results for the BA network case. The color of each grid point represents the proportion of successful activations based on 10 numerical simulations for that parameter combination. When the system is in a stable equilibrium state on

Fig. 4.1. (b) Prediction results for the ER network case. In this case, the node activity follows a narrow distribution, as shown in Fig. 4.1. (e). The actual values are very close to the predicted values. Additionally, since the topological differences

Fig. 4.1. (c) Prediction results of the method after eliminating the power-law distribution of node activity. After removing the power-law distribution of node activity, the actual values are very close to the predicted boundary. Both Fig. 4.1. (b) and (c) are the BA network, the node activity follows a power-law distribution, as shown in Fig.4.1. (d), which does not satisfy the assumptions.

between nodes in the ER network are small, the distinction between whether activation is successful under different controls is very clear.

used to verify that the power-law distribution of node activity is the reason for the decrease in prediction accuracy. In this case, the node activity follows a narrow distribution, as shown in Fig.4.1. (f).

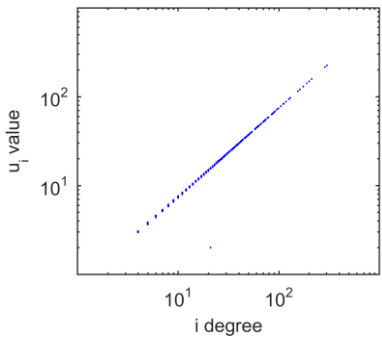

Fig.4.1. (d) Node activity follows a power-law distribution. For this model, the power-law distribution of node degree leads to a power-law distribution of node activity. Therefore, these points are concentrated along a single line.

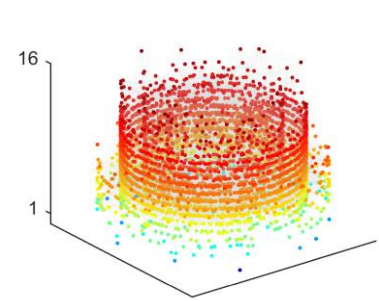

Fig.4.1. (e) In the ER network, node activity follows a narrow distribution. The positions of the nodes are determined by layering based on the shortest path. It can be observed that the values of nodes in each layer are indeed distributed within the same order of magnitude. The node values are represented by the coordinates along the z-axis.

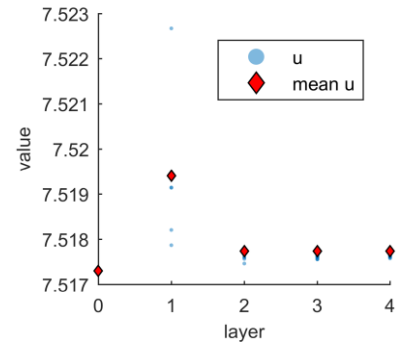

Fig.4.1. (f) The modified dynamic model, which keeps node activity as a narrow distribution in the BA network. The corresponding dynamic model is (4.2). The horizontal axis represents the layer each node belongs to, and the vertical axis shows the value of $u_i$ for each node. In this case, the activity fluctuations of each node are very small, so the values of all nodes are close to the average.

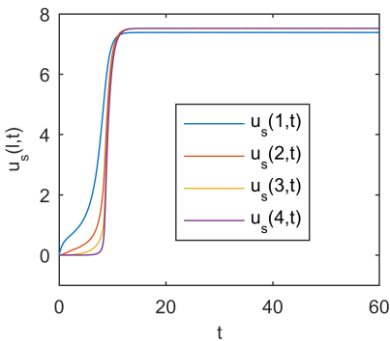

Fig.4.1. (g) The predicted system dynamics when control can activate the system. The prediction shows the change in the average state of nodes at each layer. It can be observed that the

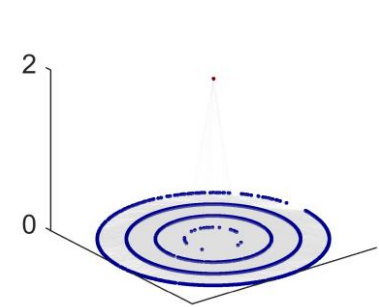

Fig.4.1. (f) The system state at the start of activation. Except for the controlled nodes, the remaining nodes are all in the low equilibrium state of 0. The nodes are layered

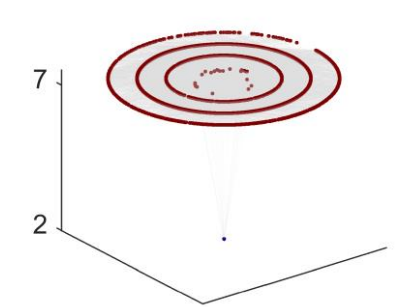

Fig.4.1. (g) The system state at the end of activation. The system state has changed from the low state in Fig.4.1. (f) to the high state. The control time is 60. From Figs.4.1. (f)

system's activation process occurs in a very short time, from a low state to a high state. The control values are $u_s = 2$ and $v_s = 2$.

Figs. 4.1. (f) and (g) display the system's changes under this control.

based on their shortest path to the selected nodes.

Figs. 4.1. (f) and (g) show the results under the BA network, after eliminating the power-law distribution of node activity. The central nodes are the controlled nodes.

and (g), it can be seen that controlling a single node has transformed the entire system from the 0 equilibrium state to the positive equilibrium state. The system activation is successful.

Figs.4.1. (a), (b), and (c) are the numerical simulation results for three different cases, aimed at verifying the explanation that node activity not following a narrow distribution leads to a decrease in the accuracy of the method. They correspond to the same predictive model. The red line indicates the boundary between the parameter ranges where activation can and cannot occur. If the method is correct, the numerical simulation results should be very close to the predicted boundary. Fig.4.1. (a) shows the numerical simulation results for the BA network. In this case, node activity follows a power-law distribution, as shown in Fig.4.1. (g). This does not meet the requirement for node activity to follow a narrow distribution. The numerical simulation results deviate somewhat from the predicted values. We believe this is due to the activity of nodes with high degree, which may not converge to the predicted values. For example, in the region close to the critical boundary where activation is predicted not to occur, nodes connected to those with average degree exhibit lower activity, and these nodes also have lower activity. However, for nodes with high degree, their interaction term is the sum of many small quantities, which does not necessarily equal a small value. Thus, the activity of high-degree nodes needs to be larger to reach equilibrium, which contradicts the prediction. To verify this explanation, we performed two additional numerical simulations. Fig.4.1. (b) shows the prediction results for an ER network, with other conditions unchanged. Since the node degrees in the ER network are similar, there is no power-law distribution of node activity. Fig.4.1. (c) shows the prediction results after modifying model (4.1) to eliminate the power-law distribution of node activity, with other conditions unchanged. The treatment is as shown in equation (4.2).

$$\frac{du_i(t)}{dt} = -B_1 u_i(t) + \frac{k}{k_i}\left(\sum_{j=i}^{N} A_{ij} \frac{v_j^2(t)}{1+v_j^2(t)}\right),$$
$$\frac{dv_i(t)}{dt} = -B_2 v_i(t) + \frac{k}{k_i}\left(\sum_{j=i}^{N} A_{ij} \frac{u_j^2(t)}{1+u_j^2(t)}\right).$$
(4.2)

Here, $k_i$ represents the degree of the $i$-th node. This approach helps prevent the degree power-law distribution of nodes from causing a power-law distribution in node activity. From this, it can be seen that when the node activity follows a narrow distribution, the prediction accuracy of the method is very high. Although the prediction accuracy decreases when the narrow distribution condition for node activity is not met, the method still has significant practical value for our purpose—activating the system. This is because the method accurately predicts whether activation is possible and the magnitude of the control needed for activation. Fig.4.1(d) shows the model (4.1) on a BA network, where node activity follows a power-law distribution. Figs.4.1(b) and (c) show narrow distributions of node activity. Fig.4.1(g) shows the changes in the average values of nodes at each layer predicted by the system. The state of each layer's nodes quickly changes from a low state to a high state during the control process. Figs.4.1(f) and (g) show the system's activation process. The controlled nodes are in the middle, the $z$-axis reflects the node's state, and the colors represent different states.

4.3 Reciprocity Model

Reciprocal relationships are a common type of interaction. They are widely present in biological and social systems. Scholars have conducted in-depth research on these relationships. For example, the mutualism between flowering plants and pollinators, and between businesses and consumers. In the following, we use the model below for numerical simulation:

$$\frac{du_i(t)}{dt} = -au_i(t) + \frac{b\sum_{j=1}^{N} A_{ij} v_j^2(t)}{1 + c\sum_{j=1}^{N} A_{ij} v_j^2(t)},$$

$$\frac{dv_i(t)}{dt} = -dv_i(t) + \frac{e\sum_{j=1}^{N} A_{ij} u_j^2(t)}{1 + f\sum_{j=1}^{N} A_{ij} u_j^2(t)},$$

(4.3)

Here, the parameters $a$, $b$, $c$, $d$, $e$, and $f$ are all positive numbers. $u_i(t)$ and $v_i(t)$ represent the values of node $i$ at time $t$. $a$ and $d$ represent the decay rates when only the nodes themselves exist. The expression for the interaction term shows that, in this system, the growth of node degree does not lead to infinite growth of the interaction term. Therefore, in a BA network where node degrees follow a power-law distribution, the prediction method is also effective. This model is inspired by the one proposed by Jiang et al. In their system, some components of the equilibrium point are zero, while others are positive. We believe this does not match many real-world reciprocal systems, as reciprocal systems usually rely on the existence of one group as a prerequisite for the existence of another. Hence, we have improved the original reciprocity model. The dynamical system that eliminates the effect of node degree on node activity is shown in equation (4.4):

$$\frac{du_i(t)}{dt} = -au_i(t) + \frac{\frac{bk}{k_i}\sum_{j=1}^{N} A_{ij}v_j^2(t)}{1 + \frac{ck}{k_i}\sum_{j=1}^{N} A_{ij}v_j^2(t)},$$

$$\frac{dv_i(t)}{dt} = -dv_i(t) + \frac{\frac{ek}{k_i}\sum_{j=1}^{N} A_{ij}u_j^2(t)}{1 + \frac{fk}{k_i}\sum_{j=1}^{N} A_{ij}u_j^2(t)}.$$

(4.4)

The parameter values are $a = 5$, $b = 4$, $c = 0.5$, $d = 3$, $e = 3$, and $f = 0.5$. The results of the numerical simulation are shown in Fig.4.2.

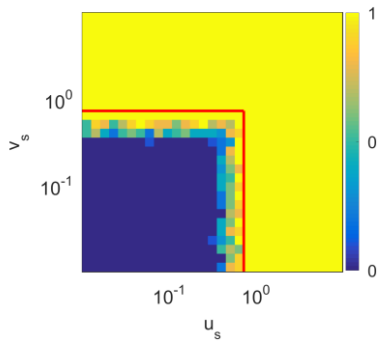

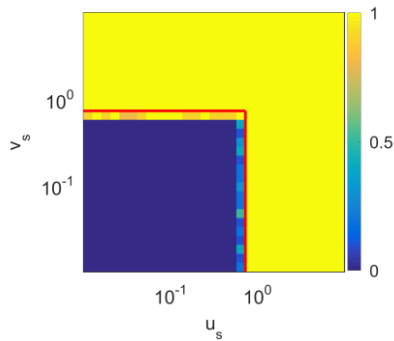

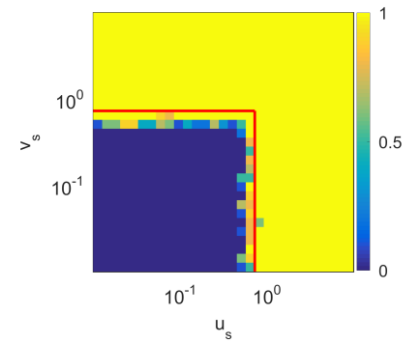

Fig.4.2(a) shows the prediction results for the BA network. The color of each grid point represents the proportion of successful activations in ten numerical simulations for that parameter combination.

Fig.4.2(b) shows the prediction results for the ER network. Due to the small topological differences between nodes in the ER network, the distinction in activation success under different controls is quite clear.

Fig.4.2(c) shows the prediction results when the power-law distribution of node activity is eliminated. The numerical simulation results are very close to the predicted boundary.

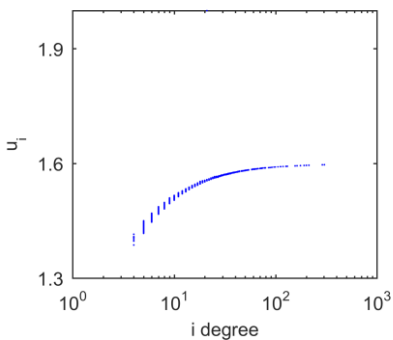

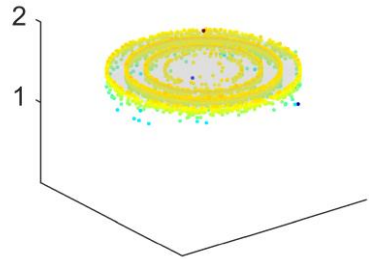

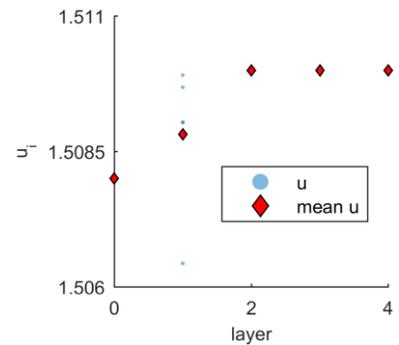

Fig.4.2(d) The distribution of node activity in the BA network. The horizontal axis represents node degree, and the vertical axis represents node activity. Node activity follows a narrow distribution. For this model, changes in node

Fig.4.2(e) In the ER network, node activity follows a narrow distribution. The positions of the nodes are determined based on the shortest path hierarchy. It can be observed that the values of nodes in each layer are

Fig.4.2(f) The node states of the modified dynamical model. The effect of node degree on activity has been eliminated. The corresponding dynamical model is (4.4). The horizontal axis represents the layer to which each node

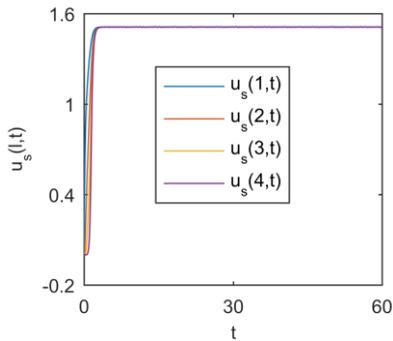
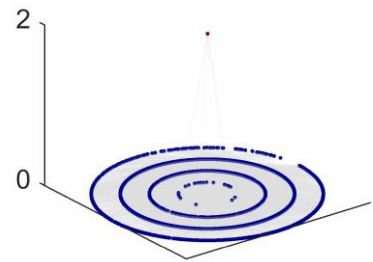
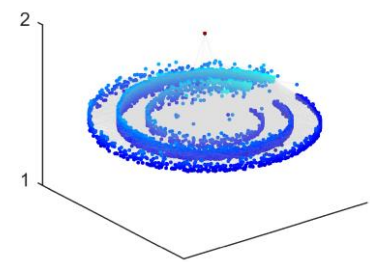

| | | |
|---|---|---|
| degree do not cause large variations in node activity. | indeed distributed within the same order of magnitude. Node values are represented by the coordinates along the $z$-axis. | belongs, and the vertical axis shows the value of $z$ for each node. |
| Fig.4.2(g) The predicted system changes when control can activate the system. The prediction is for the change in the average state of nodes at each layer. It can be seen that the system's activation process occurs within a very short time, transitioning from a low state to a high state. The controls are set to $u_s=2$ and $v_s=2$. Figs.4.2(f) and (g) display the system's changes under this control. | Fig.4.2(f) The system state at the start of activation. Except for the controlled node, all other nodes are in the low equilibrium state of 0. The nodes are layered based on their shortest path to the selected node. Figs.4.2(f) and (g) both show the numerical simulation results of model (4.3) under the BA network. The central node is the controlled node. | Fig.4.2(g) shows the system state at the end of activation. The system state has changed from the low state in Fig.4.2(f) to a high state. The control time is 60. From Figs.4.2(f) and (g), it can be seen that controlling a single node has shifted the entire system from a zero equilibrium state to a positive equilibrium state. The system activation is successful. |

Figs.4.2(a), (b), and (c) show the results of ten numerical simulations for each control combination under three different conditions. These three conditions all meet the requirements, and the corresponding prediction models are the same. The red lines represent the critical lines for activation and non-activation. From these, it can be seen that the predictions for all three cases are accurate. Figs.4.2 and 4.1 serve a similar purpose, both aiming to demonstrate the prediction performance of the method. Figs.4.2(d), (e), and (f) show that the activity of nodes is narrowly distributed in these three cases. Fig.4.2(g) presents the predicted results for the activation process. A total of 54,000 numerical simulations were carried out in Figs.4.1 and 4.2, thoroughly validating the effectiveness of the method.